\documentclass[a4paper,11pt]{article}
\usepackage{amssymb, amsmath, amsthm, amscd, amsfonts}

\usepackage[]{amsmath}
\newtheorem{thm}{Theorem}

\newtheorem{xm}{Example}
\newtheorem{al}{Algorithm}
\newenvironment{pr}[1][Proof]{\noindent\textbf{#1.} }{\ \rule{0.5em}{0.5em}}

\newtheorem{definition}{Definition}
\begin{document}
\title{\textbf{A note on \lq\lq A simple algorithm to search for all MCs in networks"}}
\author{\textbf{Majid~Forghani-elahabad\footnote{Email
address:~forghanimajid@mehr.sharif.ir} ,\
 Nezam ~Mahdavi-Amiri\footnote{Corresponding author: Fax:++982166005117, Phone: ++982166165607, Email
address: nezamm@sina.sharif.edu} }}
\date{}
\maketitle $\vspace{-14mm}$
\begin{center}
\textit{{\scriptsize Faculty of Mathematical Sciences, Sharif University of 
Technology, Tehran, Iran}}\\
\end{center}

\hrule \section*{ Abstract} $\indent$ 

Recently, Yeh [Yeh, WC. (2006). \emph{A simple algorithm to search for all MCs in networks}. European Journal of Operational Research, 174,  1694--1705.] has proposed a simple algorithm to find all the Minimal Cuts in an undirected graph. However, the algorithm does not work properly.
 Here, using an example, a defect of this algorithm is illustrated, and then the corresponding result is shown to be incorrect. Moreover, a correct version of the algorithm is established.

\vspace{2mm} \noindent {\normalsize{Keywords:}}
Reliability; Minimal Cut (MC); Algorithm\\
\hrule
\section{ Introduction}
Finding all Minimal Cuts ({\it MC}s) has extensive application in network reliability (Levitin et al., 2007; Salehi and Forghani, 2009). Scanning all the {\it MC}s is an NP-hard problem (Provan and Ball, 1983). There are several algorithms for obtainment of all the {\it MC}s in a directed graph (Tan, 2003) or an undirected graph (Fard and Lee, 1999; Yeh, 2006). To determine all the {\it MC}s, Fard and Lee (1999) considered both link and node failures and proposed an algorithm deducing all the {\it MC}s of such networks making use of networks having perfect nodes. This algorithm adopts the concept of common-cause failures and does not re-enumerate {\it MC}s for the additional supposition of node failures. Tan (2003), assuming nodes with {\it k}-out-of-{\it n} property, extended the traditional directed {\it s-t} networks and, using the definition of {\it MC} for the nodes, presented an approach for finding all the {\it MC}s for all nodes, starting with the source node {\it s} and ending with the sink node {\it t}. Yeh (2006) defined an {\it MC} using a node set (called {\it MCV}), and then proposed a simple algorithm for finding all the {\it MCV}s between two special nodes {\it s} and {\it t}. Furthermore, he generalized his algorithm for obtaining all the {\it MCV}s among all pairs of nodes. Unfortunately, the algorithm proposed by Yeh (2006) has some faults and thus may not find all the {\it MCV}s of a network. Here, we first exemplify a defect of Yeh's algorithm (Yeh, 2006), and then show the corresponding result to be incorrect. Moreover, a correct version of the algorithm is provided and its correctness is established.

In the remainder of our work, in Section 2 we provide the required definitions and demonstrate flaws in the proposed algorithm and the corresponding result using an example. In Section 3, a correct form of the algorithm is established. Section 4 gives the concluding remarks. 

\section{Finding all the MCs}
Here, some required notations are described, and then an existing algorithm (Yeh, 2006) is rewritten to explain its flaws in finding all the Minimal Cuts ({\it MC}s). Moreover, the corresponding result, Lemma 5 in (Yeh, 2006), is shown to be incorrect. 

\subsection{Problem description}

Here, the same notations, nomenclature, and assumptions used by Yeh (2006) are given. Let $G(V,E)$ be a connected network with the set of nodes $V=\{s,1,2,..,n-2,t\}$ and the set of edges $E$, where $s$ and $t$ are the source and sink nodes, respectively, and $e_{uv}\in E$ denotes an edge between nodes $u$ and $v$. For an arbitrary set of nodes $U\subseteq V$, let $\overline{U}=V-U$, $E(U)=\{e_{uv}\in E| u,v\in U\}$ be the associated edges with the set of nodes $U$,  $G(U,E(U))$ be the sub network of $G(V,E)$ including only the set of nodes $U$ and its associated edges, and $MC(U)=\{e_{uv}\in E| u\in U$ and $v\in \overline{U}\}$ be the corresponding cut. An {\it MCV} candidate is a subset of nodes whose removal will cause a disconnection of nodes {\it s} and {\it t}.

\subsection{Incorrect algorithm and lemma}

Here, the proposed algorithm and its corresponding result, Lemma 5 in (Yeh, 2006), are shown to be incorrect. 

First, we believe that some results in (Yeh, 2006) are too obvious to be stated. For intuitive understanding, consider the definition of {\it MCV} as given by Yeh (2006).
\begin{definition}
An MCV candidate in $G(V,E)$, say $U$, is an MCV when MC($U$) is a minimal cut.
\end{definition}

It is known that a minimal cut, say $C$, divides the set $V$ into two separate subsets, say $U_C$ and $\overline{U}_C$, in which $s\in U_C$ and $t\in \overline{U}_C$. 

According to Definition 1, it is vividly seen that for every {\it MCV}, say $U$, {\it MC}($U$) is a minimal cut, and conversely for every minimal cut, say $C$,  $U_C$ is an {\it MCV}. Thus, $U$ is an {\it MCV} in $G(V,E)$ if and only if {\it MC}($U$) is a minimal cut. Therefore, Theorem 1 in (Yeh, 2006) could be directly deduced from definition of {\it MCV}, and there is no further need for lemmas 1-3 and Corollary 4 in (Yeh, 2006) to conclude that theorem.

\subsubsection{Flaws in the proposed algorithm in (Yeh, 2006)}

\noindent Now, some flaws in the proposed algorithm in (Yeh, 2006) are explained. For convenience, we rewrite that algorithm as Algorithm 1 below, and then show its defects using an example.

\begin{al}\hspace{-2.7cm}$\lq$ \hspace{2.1cm} '

The proposed algorithm in (Yeh, 2006) for finding all the {\it MCV}s in a network $G(V,E)$.

{\bf Step 0.} Let $i=k=0, S=U_0=\{s\}, T=V-\{s\}, N_0=\{t\},$ and $P=\phi$.

{\bf Step 1.} If there is a node $u\hspace{-.6mm}\in\hspace{-.6mm} T-N_i$ adjacent to $S$, then $S\hspace{-.2mm}\cup \hspace{-.2mm}\{u\}$ is an MCV candidate

\hspace{15mm}and go to Step 2. Otherwise, go to Step 4.

{\bf Step 2.} If $G(T-\{u\}, E(T-\{u\}))$ is a connected network, then $S\cup \{u\}$ is an MCV 

\hspace{15mm}and\hspace{-.2mm} go\hspace{-.2mm} to \hspace{-.2mm}Step 3. Otherwise, any node set containing $S\cup \{u\}$ is not an \hspace{-.2mm}MCV.

{\bf Step 3.} Let $i=i+1, k=k+1$, $U_k=S=S\cup \{u\}$, $P=P\cup \{U_k\},$ $T=T-\{u\}$, 

\hspace{15mm}$N_i=N_{i-1}$, and go to Step 1.

{\bf Step 4.} If $i=1$, then stop. Otherwise, remove the last node, say $v$, in $S$, $i=i-1$, 

\hspace{15mm}$N_i=N_i\cup \{v\}$, $T=T\cup \{v\}$, and go to Step 1.
\end{al}
The same example as given by Yeh (2006) is considered here to show how the solutions found by $\lq$Algorithm 1' depend on the order of selecting nodes in Step 1. One can see the obtained  results through the different order of node selections in Step 1 of that algorithm, in Section 6 (An example) of (Yeh, 2006).
Moreover, there is an ambiguity in the second part of Step 2 in Algorithm 1, where $G(T-\{u\}, E(T-\{u\}))$ is not connected. In fact, in such a case it has not been determined what the algorithm should perform next. In the following example, we show that when in Step 2, $G(T-\{u\}, E(T-\{u\}))$ is not connected, the algorithm leads to a wrong result, and also the transferring to any one of the steps 1, 3, or 4 will lead to unreasonable results.
\begin{xm}
Consider Fig. 1 as a network, and find all its MCVs using $\lq$Algorithm 1'.

Solution:

\vspace{2mm}
 
Step 0. Let $i=k=0, S=U_0=\{s\}, T=V-\{s\}, N_0=\{t\},$ and $P=\phi$.

Step 1. Since $T-N_0=\{1,2,3,4\}$, node 1 is selected and transfer is made to Step 2.

Step 2. Since $G(\{2,3,4,t\}, E(\{2,3,4,t\}))$ is connected, $\{s,1\}$ is found as an MCV, 

\hspace{14mm}and transfer is made to Step 3.

Step 3. Let $i=1, k=1$, $U_1=S=\{s,1\}$, $P=\{U_1\}$, $T=\{2,3,4,t\}$, $N_1=\{t\}$, and 

\hspace{14mm}transfer is made to Step 1.

Step 1. Since $T-N_0=\{2,3,4\}$, node 3 is selected and transfer is made to Step 2.

Step 2. Since $G(\{2,4,t\}, E(\{2,4,t\}))$ is not connected, the algorithm deduces that 

\hspace{14mm}any set containing set $\{s,1,3\}$ is not an MCV.

Next, if transfer is made to Step 1, the algorithm may select node 3 again and again, and so the algorithm does not terminate. If transfer is made to Step 3, $\lq$Algorithm 1' saves the set $\{s,1,3\}$ as an MCV which is an incorrect result. If transfer is made to Step 4, since $i=1$, the algorithm stops without finding the other MCVs. Hence, in any case of transfering, the algorithm is led to a unreasonable result, and consequently it cannot find all the MCVs in Fig. 1. 

\end{xm}

There are some minor and major faults in $\lq$Algorithm 1': 

(1) First, it is seen that the algorithm does not save $U_0=\{s\}$ as an {\it MCV}. To liquidate this minor weakness, Step 0 must be refined as follows:

\vspace{2mm}

$\lq${\bf Step 0.} Let $i=k=0, S=U_0=\{s\}, T=V-\{s\}, N_0=\{t\},$ and $P=\{U_0\}$.'

\vspace{2mm}
(2) Another defect of the algorithm is its conclusion in the second part of Step 2. As seen in the above example, since $G(\{2,4,t\}, E(\{2,4,t\}))$ is not connected, the algorithm concludes that any set containing the set $\{s,1,3\}$ is not an {\it MCV}, whereas it is apparent that the sets $\{s,1,2,3\}$ and $\{s,1,2,3,4\}$ are {\it MCV} as well as containing the set $\{s,1,3\}$. 

\unitlength 1mm 
\linethickness{1pt}
\ifx\plotpoint\undefined\newsavebox{\plotpoint}\fi 
\begin{picture}(50,40)(10,60)

\put(64,90.5){\circle{6}} \put(64,65){\circle{6}}
\put(94,90.5){\circle{6}} \put(94,65){\circle{6}}
\put(37.75,79){\circle{6}} \put(120.5,79){\circle{6}}

\multiput(39.5,81.2)(.071799308,.03){302}{\line(1,0){.2}}
\multiput(40.2,77.5)(.05988024,-.033682635){349}{\line(1,0){.2}}
\multiput(96.7,65.5)(.06,.03){370}{\line(1,0){.2}}
\multiput(96.8,90.4)(.07,-.03){310}{\line(1,0){.2}}
\multiput(66.4,89)(.07,-.06){365}{\line(1,0){.2}}
\put(64,87.75){\line(0,-1){19.84}}
\put(94,87.75){\line(0,-1){19.84}}

\put(64,90.5){\makebox(0,0)[cc]{1}}
\put(64,65){\makebox(0,0)[cc]{2}}
\put(94,90.5){\makebox(0,0)[cc]{4}}
\put(94,65){\makebox(0,0)[cc]{3}}
\put(37.75,79){\makebox(0,0)[cc]{s}}
\put(120.5,79){\makebox(0,0)[cc]{t}}
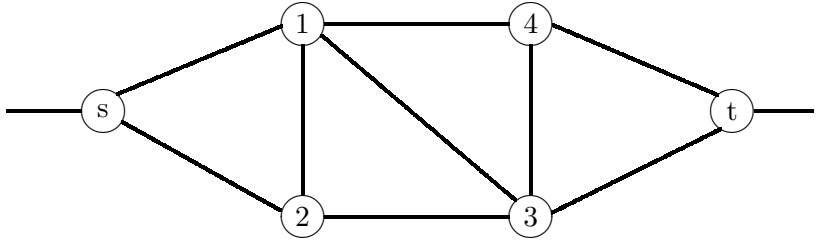
\put(80,57){\makebox(0,0)[cc]{Figure 1. An example network.}}
\put(25.1,79){\line(1,0){9.75}}
\put(123.4,79){\line(1,0){7.75}}
\put(66.85,90.5){\line(1,0){24.4}}
\put(66.85,65){\line(1,0){24.4}}

\end{picture}

\vspace{10mm}

Thus, the sentence $\lq\lq$ any node set containing $S\cup \{u\}$ is not an {\it MCV}'' must be removed from the end of Step 2.

(3) A major fault is the existent ambiguity in the second part of Step 2 in the algorithm. Likewise, in Example 1, when faced with such a state, $G(\{2,4,t\}, E(\{2,4,t\}))$ being  disconnected, transfer to any one of the steps 1,3, and 4 would be undesired. Hence, some changes are needed to rectify this shortcoming.

\subsubsection{Flaws in the proof of Lemma 5 in (Yeh, 2006)}

\noindent Yeh showed the correctness of his proposed algorithm based on use of Lemma 5 in (Yeh, 2006). Since the algorithm has some flaws, it is expected that the lemma and its proof be incorrect. Here, the proof is shown to be incorrect. First, note that the author did not define the notation $V_s(v_k)$ and used the notation $C_i$ to exhibit an {\it MCV} in the proof of the lemma, whereas the {\it MCV}s were shown by $U_i$ in the algorithm. Moreover, only the case $C_k=C_{k-1}-\{u_k\}$ was considered, while it may happen that $C_k$ is a hyper set of $C_{k-1}$. For example, in Fig. 1, $U_1=\{s,1\}$ is a hyper set of $U_0=\{s\}$. Now, the proof is shown to be incorrect. In the beginning of the proof, it is assumed that $u_k\in V_s(v_k)-C_{k-1}-T$, $k>0$, and $v_k$ is the greatest number in $C_{k-1}$ so that $V_s(v_k)-C_{k-1}-T\ \neq \phi$. The mistake is that the set $V_s(v_k)-C_{k-1}-T$ in not necessarily nonempty. For instance, take the given example in (Yeh, 2006). In the forth iteration, $U_4=\{s,1,2,3,4\}$ was found,  and in the seventh iteration, $U_5=\{s,1,2,4\}$  and $T=\{3,t\}$ were obtained. It is vividly seen that $V_s(v_5)-U_4-T$ is empty for every node $v_5\in U_4$. Therefore, the proof of Lemma 5 in (Yeh, 2006) cannot be correct. This is expected as the corresponding algorithm turns not to be correct. 

\section{The correct algorithm}
Here, using the arguments in the previous section, a correct form of the proposed algorithm in (Yeh, 2006) is presented, and then its correctness is established. According to Section 2.2.1, it is known that steps 0 and 2 and probably the stopping criterion must be refined. The assumption $P=\{U_0\}$ is used instead of $P=\phi$ in Step 0, the stated incorrect sentence is eliminated from the end of Step 2, and the stopping criterion for the algorithm is replaced by $i=0$. Moreover, to liquidate the existent ambiguity in Step 2, we use a new set $B$ and make some more changes in steps 0-2. 

\begin{al}

A correct version of the proposed algorithm in (Yeh, 2006) for finding all the MCVs in a network $G(V,E)$.\\

{\bf Step 0.} Let $i=k=0, S=U_0=\{s\}, T=V-\{s\}, N_0=\{t\},$ $B=\phi$, and $P=\{U_0\}$.

{\bf Step 1.} If there is a node $v \in T-\{B\cup N_i\}$ adjacent to $S$, then go to Step 2, else 

\hspace{15mm}go to Step 4.

{\bf Step 2.} If $G(T-\{v\}, E(T-\{v\}))$ is a connected network, let $B=\phi$ and then go to 

\hspace{15mm}Step 3, else let $B=B\cup \{v\}$ and go to Step 1.

{\bf Step 3.} Let $i\hspace{.4mm}=\hspace{.4mm}i+1, k\hspace{.4mm}=\hspace{.4mm}k+1$, $U_k\hspace{.4mm}=\hspace{.4mm}S\hspace{.4mm}=\hspace{.4mm}S\cup \{v\}$, $P\hspace{.4mm}=\hspace{.4mm}P\cup \{U_k\},$ $T\hspace{.4mm}=\hspace{.4mm}T-\{v\}$, 

\hspace{15mm}$N_i=N_{i-1}$, and go to Step 1.

{\bf Step 4.} If\hspace{.2mm} $i\hspace{-.4mm}=\hspace{-.4mm}0$, then stop, else remove the last node, $u$, in $S$, $i\hspace{-.4mm}=\hspace{-.4mm}i-1$, $N_i\hspace{-.4mm}=\hspace{-.4mm}N_i\cup \{u\}$,

\hspace{15mm}$T=T\cup \{u\}$, and go to Step 1.

\end{al}

The following theorem demonstrates the correctness of Algorithm 2.

\begin{thm}
Algorithm 2 determines all the {\it MCV}s in every connected graph.
\end{thm}
Before presenting the proof, note that, in the proof, $U_i=\{s,1,2,...,i\}$ is an {\it MCV}, $E_i=E(U_i)$ is the associated edge set with $U_i$, $\overline{U_i}=V-U_i$, $C_i$=MC($U_i$)=$\{e_{uv}|u\in U_i, v\in \overline{U_i}\}$ is the corresponding {\it MC} with $U_i$, and $V(C_i)=\{u\in U_i|e_{uv}\in C_i\}$ . It is clear that $V(C_i)\subseteq U_i$, and one can name it as the extreme nodes in $U_i$.

\vspace{2mm}

\begin{pr}
It is vividly seen that Algorithm 2 finds $U_0=\{s\}$ as the {\it MCV} in Step 0, in the beginning. Now, to the contrary, assume that there are some {\it MCV}s missed by the algorithm. Without loss of generality, assume that $U_i=\{s,1,2,...,i\}$ is a missed {\it MCV} with minimum number of nodes (note that one can rename the nodes of graph to exhibit such {\it MCV} as $U_i=\{s,1,2,...,i\}$). Since $U_i$ is an {\it MCV}, graph $G(U_i,E_i)$ is connected, and consequently there is at least one node $ j\in V(C_i)$ so that $G(U_i-\{j\},E(U_i-\{j\}))$ is connected. Since $j\in V(C_i)$, there is at least one $j'\in \overline{U_i}$ so that $e_{jj'}\in C_i$. Hence, $G(\overline{U_i}\cup\{j\}, E(\overline{U_i}\cup\{j\}))$ is connected. As a result, since $G(U_i-\{j\},E(U_i-\{j\}))$ and $G(\overline{U_i}\cup\{j\}, E(\overline{U_i}\cup\{j\}))$ are connected, $U'_i=U_i-\{j\}$ is an {\it MCV}. Because $U_i$ is a missed {\it MCV} by the algorithm with the minimum number of nodes, $U'_i$ is found by the algorithm. It is observed that Algorithm 2 finds {\it MCV}s in Step 3, and then goes to Step 1 for finding the next possible one. Therefore, if Algorithm 2 finds $U'_i$ in Step 3, then it will determine $U_i$ in the subsequent iterations by selecting $j$ as a node adjacent to $S=U'_i$. This contradicts the earlier assumption that $U_i$ is a missed {\it MCV}, and thus completes the proof.
\end{pr}

\vspace{2mm}

Here, we find all the {\it MCV}s in Fig. 1 by Algorithm 2 showing how the algorithm finds all the {\it MCV}s without missing any one of them.

\vspace{2mm}

Solution:\\

Step 0$^a$. Let $i=k=0, S=U_0=\{s\}, T=V-\{s\}, N_0=\{t\},$ $B=\phi$, and $P=\{U_0\}$.

Step 1. Since $T-\{B\cup N_0\}=\{1,2,3,4\}$, node 1 is selected, and transfer is made to 

\hspace{14mm}Step 2.

Step 2. Since $G(\{2,3,4,t\}, E(\{2,3,4,t\}))$ is connected, $B=\phi$ and transfer is made

\hspace{14mm}to Step 3.

Step 3. Let $i=1, k=1$, $U_1=S=\{s,1\}$, $P=\{U_0,U_1\}$, $T=\{2,3,4,t\}$, $N_1=\{t\}$, 

\hspace{14mm}and transfer is made to Step 1.

Step 1. Since $T-\{B\cup N_1\}=\{2,3,4\}$, node 3 is selected, and transfer is made to 

\hspace{14mm}Step 2.

Step 2$^b$. Since $G(\{2,4,t\}, E(\{2,4,t\}))$ is disconnected, $B=\{3\}$ and transfer is made

\hspace{14mm}to Step 1.

Step 1. Since $T-\{B\cup N_1\}=\{2,4\}$, node 2 is selected, and transfer is made to Step 

\hspace{14mm}2.

Step 2. Since $G(\{3,4,t\}, E(\{3,4,t\}))$ is connected, $B=\phi$ and teransfer is made to

\hspace{14mm}Step 3.

Step 3. Let $i=2, k=2$, $U_2=S=\{s,1,2\}$, $P\hspace{-.5mm}=\hspace{-.5mm}\{U_0,U_1,U_2\}$, $T=\{3,4,t\}$, $N_2=\{t\}$,

\hspace{14mm}and transfer is made to Step 1.

Step 1. Since $T-\{B\cup N_2\}=\{3,4\}$, node 3 is selected, and transfer is made to Step 

\hspace{14mm}2.

Step 2. Since $G(\{4,t\}, E(\{4,t\}))$ is connected, $B=\phi$ and transfer is made to Step 3.

Step 3. Let $i=3, k=3$, $U_3=S=\{s,1,2,3\}$, $P=\{U_0,U_1,U_2,U_3\}$, $T=\{4,t\}$, 

\hspace{14mm}$N_3=\{t\}$, and transfer is made to Step 1.

Step 1. Since $T-\{B\cup N_3\}=\{4\}$, node 4 is selected, and transfer is made to Step 2.

Step 2. Since $G(\{t\}, E(\{t\}))$ is connected, $B=\phi$ and transfer is made to Step 3.

Step 3. Let $i=4, k=4$, $U_4=S=\{s,1,2,3,4\}$, $P=\{U_0,U_1,U_2,U_3,U_4\}$, $T=\{t\}$, 

\hspace{14mm}$N_4=\{t\}$, and transfer is made to Step 1.

Step 1. Since $T-\{B\cup N_4\}=\phi$, transfer is made to Step 4.

Step 4. Since $i=4\neq 0$, let $S=S-\{4\}=\{s,1,2,3\}$, $i=i-1=3$, $N_3=\{4,t\}$, 

\hspace{14mm}$T=\{4,t\}$, and transfer is made to Step 1.

Step 1. Since $T-\{B\cup N_3\}=\phi$, transfer is made to Step 4.

Step 4. Since $i=3\neq 0$, let $S=S-\{3\}=\{s,1,2\}$, $i=i-1=2$, $N_2=N_2\cup\{3\}=$

\hspace{14mm}$\{3,t\}$, $T=\{3,4,t\}$, and transfer is made to Step 1.

Step 1. Since $T-\{B\cup N_2\}=\{4\}$, node 4 is selected, and transfer is made to Step 2.

Step 2. Since $G(\{3,t\}, E(\{3,t\}))$ is connected, $B=\phi$ and transfer is made to Step 3.

Step 3. Let $i=3, k=5$, $U_5=S=\{s,1,2,4\}$, $P=\{U_0,U_1,U_2,U_3,U_4,U_5\}$, $T=$

\hspace{14mm}$\{3,t\}$, $N_3=\{3,t\}$, and transfer is made to Step 1.

Step 1. Since $T-\{B\cup N_3\}=\phi$,  transfer is made to Step 4.

Step 4. Since $i=3\neq 0$, let $S=S-\{4\}=\{s,1,2\}$, $i=i-1=2$, $N_2=N_2\cup\{4\}=$

\hspace{14mm}$\{3,4,t\}$, $T=\{3,4,t\}$, and transfer is made to Step 1.

Step 1. Since $T-\{B\cup N_2\}=\phi$,  transfer is made to Step 4.

Step 4. Since $i=2\neq 0$, let $S=S-\{2\}=\{s,1\}$, $i=i-1=1$, $N_1=N_1\cup\{2\}=$

\hspace{14mm}$\{2,t\}$, $T=\{2,3,4,t\}$, and transfer is made to Step 1.

Step 1. Since $T-\{B\cup N_1\}=\{3,4\}$, node 3 is selected, and transfer is made to Step 

\hspace{14mm}2.

Step 2$^b$. Since $G(\{2,4,t\}, E(\{2,4,t\}))$ is disconnected, $B=\{3\}$ and transfer is made

\hspace{14mm}to Step 1.

Step 1. Since $T-\{B\cup N_1\}=\{4\}$, node 4 is selected, and transfer is made to Step 2.

Step 2. Since $G(\{2,3,t\}, E(\{2,3,t\}))$ is connected, $B=\phi$ and teransfer is made to

\hspace{14mm}Step 3.

Step 3. Let $i=2, k=6$, $U_6=S=\{s,1,4\}$, $P=\{U_0,U_1,U_2,U_3,U_4,U_5,U_6\}$, 

\hspace{14mm}$T=\{2,3,t\}$, $N_2=\{2,t\}$, and transfer is made to Step 1.

Step 1. Since $T-\{B\cup N_2\}=\{3\}$, node 3 is selected, and transfer is made to Step 2.

Step 2$^b$. Since $G(\{2,t\}, E(\{2,t\}))$ is disconnected, $B=\{3\}$ and transfer is made to

\hspace{14mm}Step 1.

Step 1. Since $T-\{B\cup N_2\}=\phi$, transfer is made to Step 4.

Step 4. Since $i=2\neq 0$, let $S=S-\{4\}=\{s,1\}$, $i=i-1=1$, $N_1=N_1\cup\{4\}=$

\hspace{14mm}$\{2,4,t\}$, $T=\{2,3,4,t\}$, and transfer is made to Step 1.

Step 1. Since $T-\{B\cup N_1\}=\phi$, transfer is made to Step 4.

Step 4$^c$. Since $i=1\neq 0$, let $S=S-\{1\}=\{s\}$, $i=i-1=0$, $N_0=N_0\cup\{1\}=$

\hspace{14mm}$\{1,t\}$, $T=\{1,2,3,4,t\}$, and transfer is made to Step 1.

Step 1. Since $T-\{B\cup N_0\}=\{2,4\}$, node 2 is selected, and transfer is made to 

\hspace{14mm}Step 2.

Step 2. Since $G(\{1,3,4,t\}, E(\{1,3,4,t\}))$ is connected, $B=\phi$ and transfer is made

\hspace{14mm}to Step 3.

Step 3. Let $i=1, k=7$, $U_7=S=\{s,2\}$, $P=\{U_0,U_1,U_2,U_3,U_4,U_5,U_6,U_7\}$, 

\hspace{14mm}$T=\{1,3,4,t\}$, $N_1=\{1,t\}$, and transfer is made to Step 1.

Step 1. Since $T-\{B\cup N_1\}=\{3,4\}$, node 3 is selected, and transfer is made to Step 

\hspace{14mm}2.

Step 2. Since $G(\{1,4,t\}, E(\{1,4,t\}))$ is connected, $B=\phi$ and transfer is made to 

\hspace{14mm}Step 3.

Step 3. Let $i=2, k=8$, $U_8=S=\{s,2,3\}$, $P=\{U_0,U_1,U_2,U_3,U_4,U_5,U_6,U_7,U_8\}$, 

\hspace{14mm}$T=\{1,4,t\}$, $N_2=\{1,t\}$, and transfer is made to Step 1.

Step 1. Since $T-\{B\cup N_2\}=\{4\}$, node 4 is selected, and transfer is made to Step 2.

Step 2$^b$. Since $G(\{1,t\}, E(\{1,t\}))$ is disconnected, $B=\{4\}$ and transfer is made to

\hspace{14mm}Step 1.

Step 1. Since $T-\{B\cup N_2\}=\phi$, transfer is made to Step 4.

Step 4. Since $i=2\neq 0$, let $S=S-\{3\}=\{s,2\}$, $i=i-1=1$, $N_1=N_1\cup\{3\}=$

\hspace{14mm}$\{1,3,t\}$, $T=\{1,3,4,t\}$, and transfer is made to Step 1.

Step 1. Since $T-\{B\cup N_2\}=\phi$, transfer is made to Step 4.

Step 4. Since $i=1\neq 0$, let $S=S-\{2\}=\{s\}$, $i=i-1=0$, $N_0=N_0\cup\{2\}=$

\hspace{14mm}$\{1,2,t\}$, $T=\{1,2,3,4,t\}$, and transfer is made to Step 1.

Step 1. Since $T-\{B\cup N_0\}=\{3\}$, there is no node adjacent to $S$ and transfer is

\hspace{14mm}made to Step 4.

Step 4. Since $i=0$, the algorithm stops.\\

$\bullet^a$: Note that Algorithm 2 found $U_0=\{s\}$ owing to the considered change in Step 0 of the algorithm.

$\bullet^b$: The usage of the new set $B$ and modification of Step 2 rectified the mentioned defect of the algorithm proposed by Yeh (2006), and Algorithm 2 acted correctly when the sub network $G(T-\{v\}, E(T-\{v\}))$ was disconnected.

$\bullet^c$: The change of stopping criterion of the algorithm caused the algorithm not to stop before finding all the {\it MCV}s.

\section*{Conclusions}

There are a number of algorithms for finding all the Minimal Cuts ({\it MC}s) considering different categories of networks such as directed or undirected networks, network with node failures, link failures, or both node and link failures, networks with {\it k}-out-of-{\it n} nodes, cyclic or acyclic networks, etc. Here, an algorithm proposed by Yeh [Yeh, WC. (2006). \emph{A simple algorithm to search for all MCs in networks}. European Journal of Operational Research, 174,  1694--1705] was shown to be incorrect, demonstrating flaws in the correctness proof. A modified version of the algorithm was proposed and its correctness was established. 

\subsection*{Acknowledgements}
The authors thank the Research Council of Sharif University of Technology for its support.

{}
\end{document}